\documentclass{article} 
\usepackage{amssymb, amsfonts, amsmath, amsthm}
\usepackage{epsfig,lpic,wrapfig}
\usepackage[T2A,T1]{fontenc}
\usepackage[utf8]{inputenc}
\usepackage[english,russian]{babel}
\usepackage{graphicx}
\usepackage{url}
\usepackage{tikz}
\usetikzlibrary{automata,chains}
\usepackage[top=0.9in, bottom=0.9in,left=0.9in, right=0.9in, paperwidth=6in, paperheight=9in]{geometry}
\usepackage[colorlinks=true,
citecolor=black,
linkcolor=black,
anchorcolor=black,
filecolor=black,
menucolor=black,
urlcolor=black%
]{hyperref}
\hypersetup{pdftitle={Блужданья по цепям},
pdfauthor={Александр Гиль и Антон Петрунин}}
\def\man{\includegraphics[scale=.021]{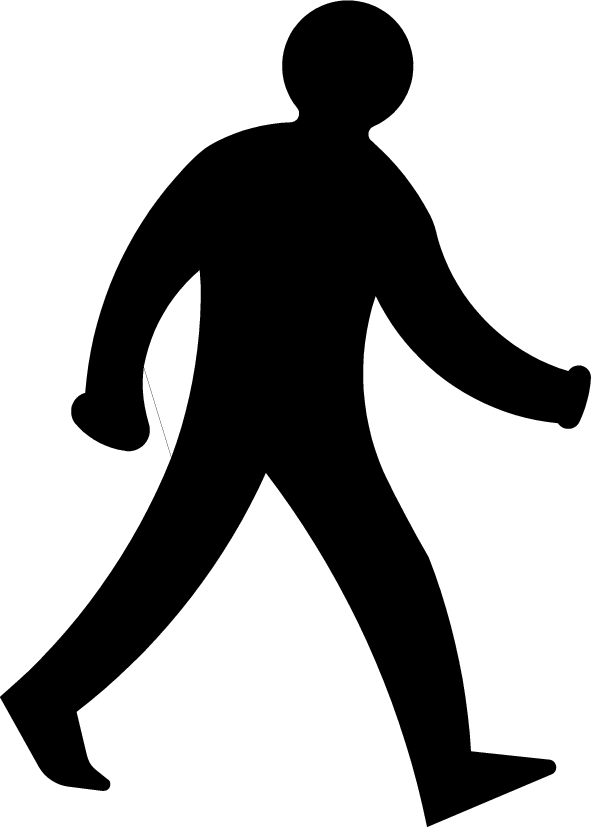}}
\def\dog{\includegraphics[scale=.020]{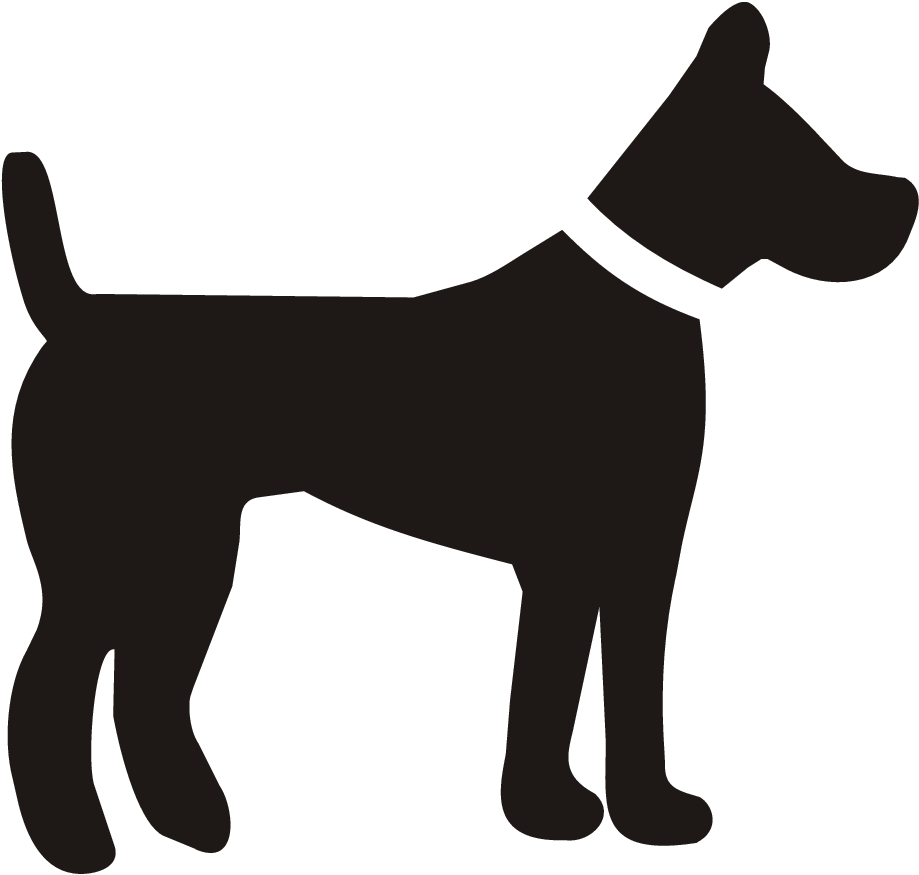}}
\def\home{\includegraphics[scale=.021]{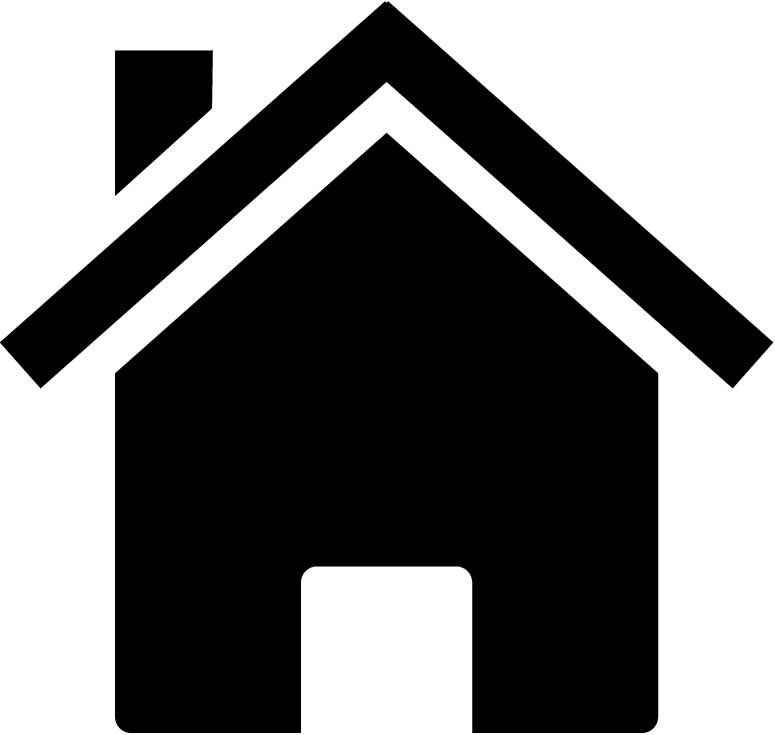}}
\def\beer{\includegraphics[scale=.008]{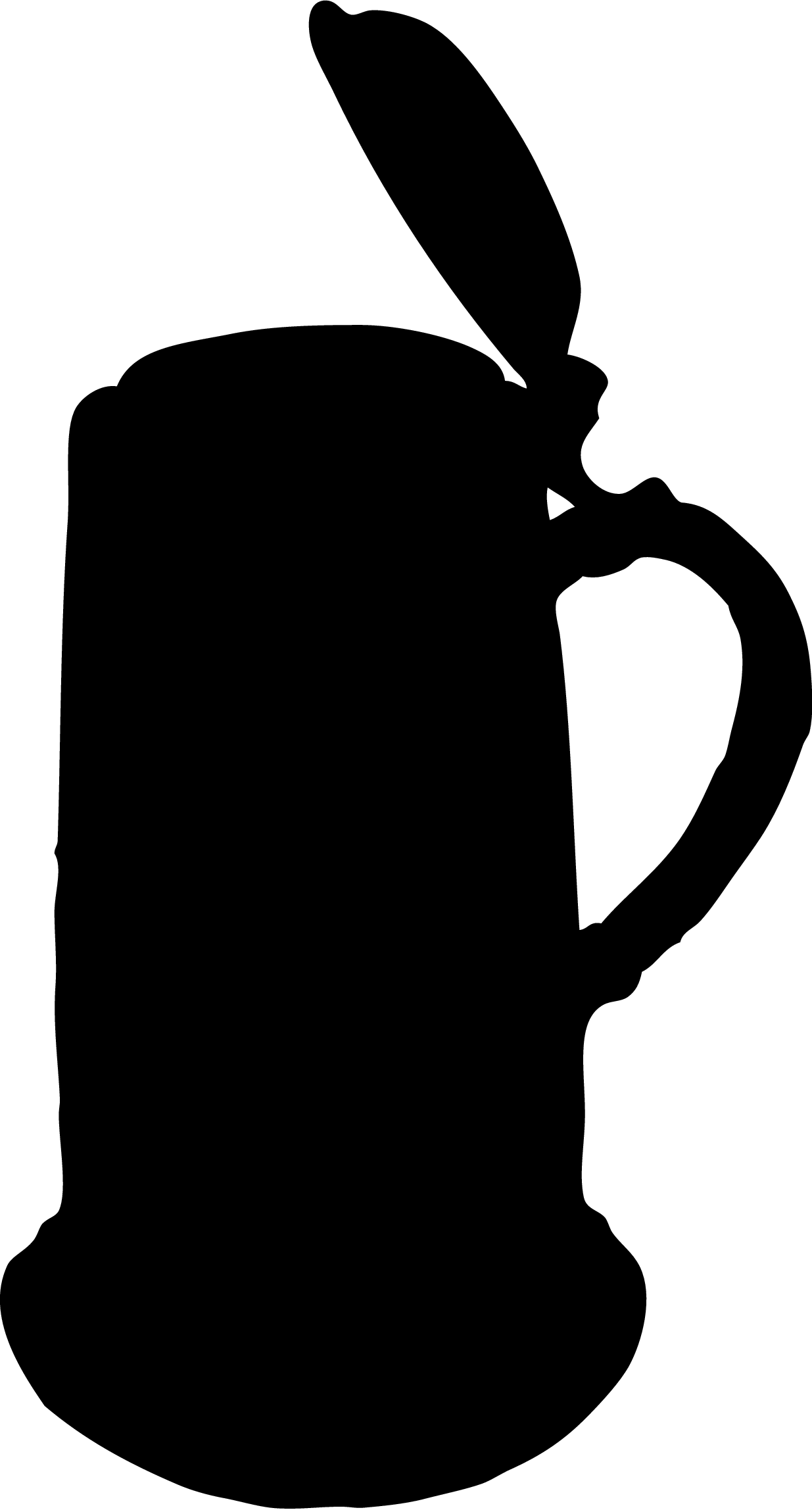}}
\def\robot{\includegraphics[scale=.1]{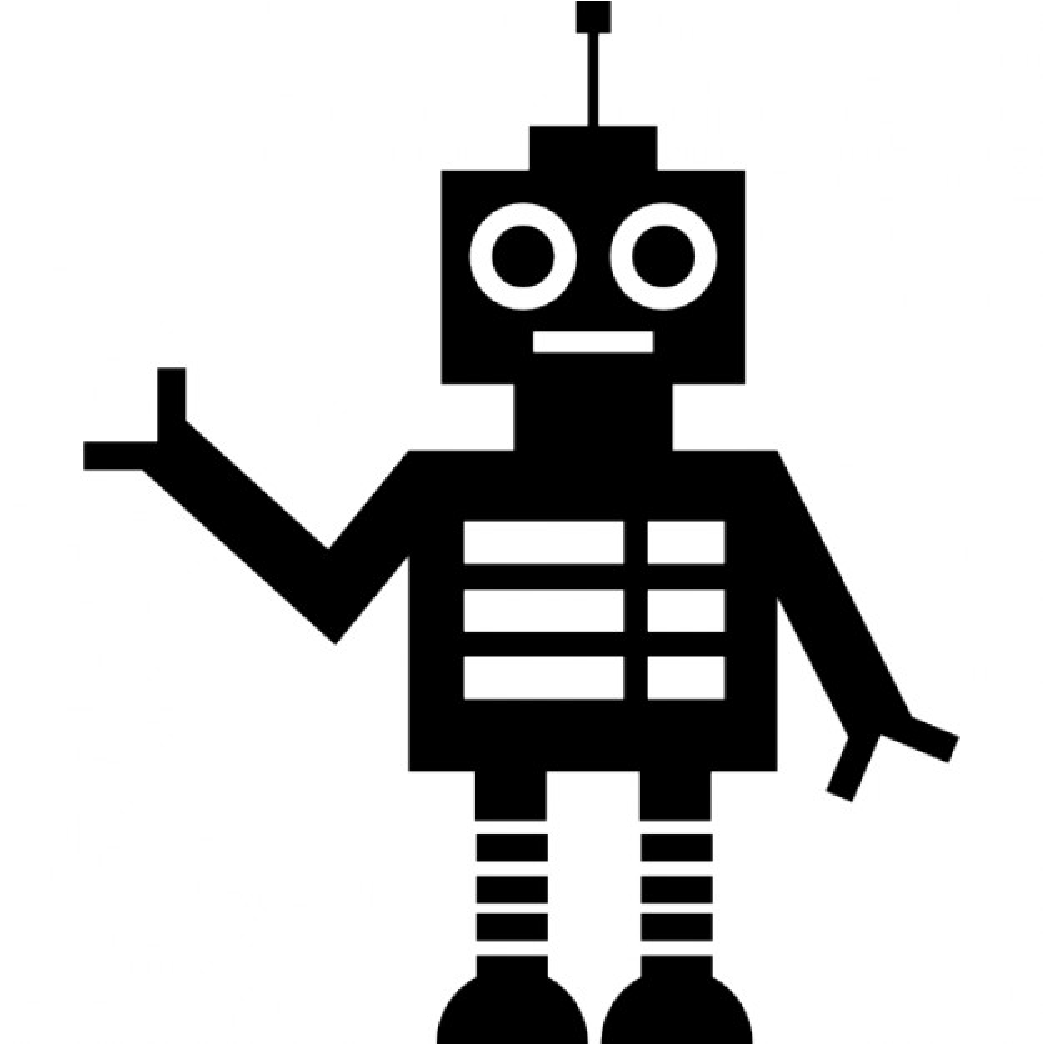}}
\usepackage{indentfirst}

\newcommand{\nofootnote}[1]{%
\begingroup\def\thefootnote{}\footnotetext{#1}\endgroup}
\begin{document}
\title{Блужданья по цепям}
\author{Александр Гиль и Антон Петрунин}
\date{}
\maketitle

\nofootnote{Заметка основана на лекции, прочитанной первым автором в школьном кружке «Northwest Academy of Sciences» в Сиэтле.
Авторы признательны 
Арсению Акопяну, 
Андрею Бураго,
Александру Ващилло,
Алексею Новикову 
и Сергею Табачникову за помощь.}
\section{Вероятность и математическое ожидание}

Подкинув игральный кубик, мы можем выбросить с равными шансами 1, 2, 3, 4, 5 или 6 очков.
\emph{Исход} (результат) такого \emph{испытания} (число очков, выпавшее на 
верхней грани кубика) является простым примером \emph{случайной величины}.
Естественно предположить, что исход одного испытания (подкидывания кубика)
не зависит от исходов других таких же испытаний.

Давайте записывать значения нашей случайной величины, много раз подкидывая кубик,
и вычислять долю испытаний, давших конкретный исход, например, 5 очков.
При неограниченном увеличении числа испытаний эта доля стремится к числу, 
которое называется \emph{вероятностью} этого исхода.
Поскольку шансы каждого из шести исходов одинаковы, вероятность каждого исхода равна $\tfrac16$ --- сумма долей шести возможных исходов всегда равна 1, и в пределе все эти шесть долей уравниваются.

Давайте теперь не только записывать результат каждого испытания, но и считать среднее
арифметическое всех записанных на данный момент результатов. Если эта последовательность стремится к
определённому числу, это число называется \emph{математическим ожиданием} или \emph{средним значением}
случайной величины (или просто \emph{средним}).

В нашем примере искомое среднее значение 
числа очков можно посчитать по формуле
\[\tfrac16\cdot1+\tfrac16\cdot2+\tfrac16\cdot3+\tfrac16\cdot4+\tfrac16\cdot5+\tfrac16\cdot6=3\tfrac12.\]
Действительно, в пределе на долю каждого исхода 1, 2, 3, 4, 5 и 6
приходится $\tfrac16$ количества всех исходов, и, значит, их среднее арифметическое должно 
стремиться к левой стороне равенства.

Иначе говоря, для вычисления среднего значения мы должны вычислить \emph{взвешенную сумму} значений нашей величины
для каждого исхода (1, 2, 3, 4, 5 и 6), взяв вероятность каждого исхода ($\tfrac16$) в качестве его веса.

\medskip

Хотя изложенные выше утверждения и кажутся очевидными, для того,
чтобы придать им точный математический смысл,
требуется развить теорию, которая останется за рамками нашей заметки.
Мы обсудим только способы нахождения вероятностей и средних значений в чуть более сложных ситуациях.

\medskip

Вообще говоря, может оказаться, что среднее значение (математическое ожидание), вычисляемое как предел для неограниченно продлеваемой серии испытаний, не существует для данной схемы испытаний. 
Такое может случиться,
поскольку у бесконечной последовательности чисел предельного значения может и не существовать.

Знаменитым примером ситуации без предельного значения является «Санкт-Петербургская лотерея», состоящая в следующем: банкомёт подбрасывает  монету несколько раз, пока не выпадает орёл.
Если орёл выпадает на первом броске, банкомёт выплачивают игроку 1 дукат;
если орёл выпадает при втором броске, то 2 дуката;
если при третьем --- 4 дуката, и так далее (на каждом новом броске ставка удваивается). 
Средний выигрыш в такой лотерее оказывается бесконечным.
Об этом парадоксе и его истории можно прочитать в статье \cite{kudryavtzev}.
В наших задачах такого происходить не будет,
но мы не будем это строго доказывать.

\section{Робот Чебуратор на коротком столе} 

Электромеханическая игрушка «Робот Чебуратор», 
будучи включённой, через равные интервалы времени делает шаги одинаковой длины налево или направо,
случайно выбирая между этими двумя направлениями с одинаковой вероятностью $\tfrac12$.

{

\begin{wrapfigure}{r}{19 mm}
\begin{lpic}[t(10 mm),b(0 mm),r(0 mm),l(0 mm)]{pics/stol(1)}
\lbl[b]{6,4;\robot}
\end{lpic}
\end{wrapfigure}

Предположим, Чебуратор стоит на левом краю короткого стола: шагая влево, он чебурахнется со стола.
Однако, он может шагнуть вправо и остаться на столе (на правом краю) --- но чебурахнется после второго шага направо. 
Нас интересуют следующие две задачи.

}
\begin{itemize}
\item Зная начальное положение Чебуратора, 
какова вероятность того, что робот чебурахнется с левого края стола, и какова --- что с правого?
\item Сколько в среднем шагов сделает робот, чебурахнувшись на последнем шаге?
\end{itemize}
При этом мы полагаем, что существование требуемых чисел в обоих вопросах обосновано
(это действительно верно, хотя требует доказательства); нам остаётся лишь найти значения этих чисел.

\medskip
\noindent\textit{Решение первой задачи.}
Пусть $p_1$ и $p_2$ обозначают вероятности чебураханья с правого края стола из первого (на левом краю стола) и второго (на правом краю стола) положения Чебуратора.
Заметим, что
\[p_1=\tfrac12\cdot0+\tfrac12\cdot p_2.\]
Действительно, на первом шаге с вероятностью $\tfrac12$ робот чебурахнется налево,
и, значит, шансов чебурахнуться направо у него уже не будет; 
отсюда слагаемое $0=\tfrac12\cdot0$. 
С той же вероятностью $\tfrac12$ он перейдёт на второе место (на правом краю стола), где вероятность чебурахнутся направо мы обозначили $p_2$; отсюда слагаемое $\tfrac12\cdot p_2$. 

То есть, зная две вероятности чебураханья направо после двух вариантов первого шага (0 и $p_2$), мы посчитали их взвешенную сумму, взяв в качестве весов вероятности соответствующих вариантов первого шага ($\tfrac12$ и $\tfrac12$).

Аналогично, рассмотрев два возможных исхода шага робота из второй позиции (на правом краю стола), получаем уравнение 
\[p_2=\tfrac12\cdot p_1+\tfrac12\cdot 1.\]

Систему из этих двух уравнений с двумя неизвестными ($p_1$ и $p_2$) легко решить алгебраически, но мы сделаем это по-другому (этот другой способ пригодится нам позже при решении более сложной задачи). 

Первое уравнение говорит нам, что $p_1$ равно среднему арифметическому 0 и $p_2$; второе уравнение говорит нам, что $p_2$ равно среднему арифметическому $p_1$ и 1. То есть, в последовательности из четырёх чисел
\[0,\  p_1,\  p_2,\ 1\]
каждое из чисел, кроме 0 и 1 на краях, является средним арифметическим своих соседей, а, значит, эта последовательность является арифметической прогрессией. 
Интервал от 0 до 1 разбит на три равных промежутка между членами последовательности; значит, $p_1=\tfrac13$ и $p_2=\tfrac23$.

Заметим, что если $q_1$ и $q_2$ --- вероятности чебураханья налево (при первом и втором начальных положениях), 
то, проведя аналогичные вычисления, мы получим $q_1=\tfrac23$ и $q_2=\tfrac13$.
Этот же результат можно получить, посмотрев на задачу через зеркало, меняющее местами право и лево.

Таким образом, из первой (левой) позиции робот чебурахается направо и налево с вероятностями $\tfrac13$ и $\tfrac23$.
\qed
\medskip

Заметим, что 
\[p_1+q_1=p_2+q_2=1,\]
а, значит, вероятность того, что робот никогда не чебурахнется (строго чередуя шаги направо и налево), равна нулю. 
К этому же выводу мы бы пришли, вычислив вероятность того, что Чебуратор остаётся на столе после шага номер $n$, напрямую: эта вероятность равна $(\tfrac12)^n$ и приближается сколь угодно близко к нулю с возрастанием $n$ --- то есть, равна нулю в пределе при неограниченном росте $n$.
\medskip

Аналогично решается и вторая задача.

\medskip
\noindent\textit{Решение второй задачи.}
Обозначим $s_1$ и $s_2$ средние количества шагов, которые сделает робот, чебурахнувшись на последнем шаге,
начиная с первой или второй позиции соответственно.
Взглянув на робота через зеркало, меняющее местами право и лево, мы убеждаемся, что $s_1=s_2$ ---
но мы и так это скоро увидим из вычислений.

Заметим, что 
\[s_1=\tfrac12\cdot1+\tfrac12\cdot (1+s_2).\]
Действительно, на первом шаге с первой позиции
с вероятностью $\tfrac12$ робот чебурахнется налево, то есть он проделает всего один шаг;
отсюда слагаемое $\tfrac12\cdot1$. С той же вероятностью он шагнёт вправо на вторую позицию, откуда в среднем он пройдёт ещё $s_2$ шагов. Учтя начальный шаг, получаем второе слагаемое
$\tfrac12\cdot (1+s_2)$.

Аналогично получаем равенство
\[s_2=\tfrac12\cdot (1+s_1)+\tfrac12\cdot 1.\]
Мы опять получили систему из двух уравнений с двумя неизвестными $s_1$ и $s_2$.
Решив эту систему, получаем 
$s_1=s_2=2$.
\qed
\medskip

Попробуйте решить следующую задачу тем же способом.

\begin{wrapfigure}{r}{35 mm}
\begin{lpic}[t(-0 mm),b(-2 mm),r(0 mm),l(0 mm)]{pics/gorod(1)}
\lbl[br]{6,8;\home}
\lbl[b]{12,17;\man}
\lbl[t]{27.5,5.5;\dog}
\end{lpic}
\end{wrapfigure}

\medskip
\noindent\textbf{Задача.}
На рисунке вы видите кусок, вырезанный из плана города, с отмеченными на нём собакой (\dog), 
её забывчивым хозяином (\man) и их домом (\home).

Хозяин ищет свой дом, 
выбирая с равной вероятности одну из сходящихся на перекрестке дорог, 
и идёт до следующего перекрёстка ровно одну минуту.
Дойдя до следующего перекрёстка, он опять выбирает с равной вероятностью 
одну из дорог (возможно, ту же, по которой только что шёл), 
и идёт по ней --- и так далее, пока не попадёт на тот перекрёсток, на котором его дом.

Собака двигается по тому же принципу, но бегает в полтора раза быстрее хозяина.

Докажите, что из данного начального положения собака в среднем прибежит домой на 20 секунд раньше хозяина, независимо от устройства остальной части города.

\section{Решения прямым подсчётом}

Две задачи про Чебуратора, рассмотренные выше, допускают решения прямым подсчётом вероятностей. 
Мы опишем их, используя те же обозначения, что и раньше.

\medskip
\noindent\textit{Решения.}
Решая вторую задачу, заметим, что робот чебурахается на $n$-ом шаге либо влево, либо вправо с вероятностью $(\tfrac12)^n$.
Значит, среднее значение числа шагов можно найти, просуммировав ряд
\[s_1=s_2=\tfrac12\cdot1+(\tfrac12)^2\cdot 2+(\tfrac12)^3\cdot 3+\dots.\]
Для суммирования такого ряда требуется некоторое умение; однако, ряд действительно сходится к сумме, равной $2$ --- и неудивительно, ведь это же значение было получено выше.

Далее заметим, что если робот начинает движение с левой позиции, 
то он может чебурахнуться направо только на чётных шагах --- 
а если он начинает с правой позиции, то только на нечётных.
То есть,
\begin{align*}
p_1=(\tfrac12)^2+(\tfrac12)^4+(\tfrac12)^6+\dots
\\
p_2=(\tfrac12)^1+(\tfrac12)^3+(\tfrac12)^5+\dots
\end{align*}
Применив формулу для суммы геометрической прогрессии, получаем те же результаты: $p_1=\tfrac13$ и $p_2=\tfrac23$.
\qed
\medskip

Как видите, этот способ оказался сложней.
Кроме того, он не допускает лёгкого обобщения на случай длинного стола (на котором Чебуратор может сделать более одного шага).
Как мы увидим ниже, в этом случае наше первое решение остаётся практически без изменений.

То, что среднее значение числа шагов с чебураханьем на последнем шаге
равно двум, можно также увидеть, поставив следующий мысленный эксперимент.
Представьте себе, что при каждом чебураханьи мы ставим робота обратно на стол.
Тогда вероятность того, что робот чебурахнется на шаге $n$,
равна $\tfrac12$,
то есть робот будет чебурахаться в среднем на каждом втором шаге.
Значит, среднее количество шагов между чебураханьями равно 2.

\section{Робот Чебуратор на длинном столе} 

Представим теперь, что Чебуратор стоит на длинном столе.

\begin{center}
\begin{lpic}[t(15 mm),b(0 mm),r(2 mm),l(1 mm)]{pics/dlinnyj-stol(1)}
\lbl[b]{14.75,4;\robot}
\lbl[b]{-.5,4;$\curvearrowleft$}
\lbl[b]{5.5,4;$\curvearrowleft$}
\lbl[b]{11.5,4;$\curvearrowleft$}
\lbl[b]{20.5,4;$\curvearrowright$}
\lbl[b]{26.5,4;$\curvearrowright$}
\lbl[b]{32.5,4;$\curvearrowright$}
\lbl[b]{38.5,4;$\curvearrowright$}
\end{lpic}
\end{center}

Предположим, что робот чебурахнется с левого края стола, сделав $i$ шагов налево,
и с правого края стола, сделав $j$ шагов направо.

Нас интересуют те же задачи:

\begin{itemize}
\item Какова вероятность того, что робот чебурахнется с левого края стола, и какова --- что с правого?
\item Сколько шагов в среднем сделает робот, чебурахнувшись на последнем шаге?
\end{itemize}

\medskip
\noindent\textit{Решение первой задачи.}
Заметим, что сумма $m=i+j$ --- одна и та же для всех позиций, и число позиций равно $m-1$.
Пронумеруем возможные положения Чебуратора числами от $1$ до $m-1$
по числу шагов налево, которых Чебуратору нужно сделать, чтобы чебурахнуться.

Нам будет удобно добавить ещё два положения под номерами $0$ и $m$:
первое соответствует тому, что робот чебурахнулся налево, 
а второе соответствует тому, что робот чебурахнулся направо
(из этих позиций выхода нет).

Пусть $p_i$ есть вероятность того, что Чебуратор чебурахнется направо,
если начинает с позиции под номером $i$.
Естественно, мы имеем 
\[p_0=0\ \  \text{и}\ \  p_{m}=1.\]

Применим тот же метод, что и раньше.
Стоя на позиции номер $i$, с вероятностью $\tfrac12$ Чебуратор перейдёт на позицию номер $i+1$.
В этом случае вероятность того, что он чебурахнется вправо, будет $p_{i+1}$.
С той же вероятностью $\tfrac12$ он перейдёт на позицию номер $i-1$. 
В последнем случае вероятность того, что он чебурахнется вправо, будет $p_{i-1}$.
То есть,
\[p_i=\tfrac12\cdot p_{i-1}+\tfrac12\cdot p_{i+1}.\]

Иначе говоря, в последовательности из $m+1$ числа
\[0=p_0,\ p_1,\ p_2,\dots,\ p_{m}=1\] 
каждое число, кроме двух крайних, является средним арифметическим 
соседей; мы получили арифметическую прогрессию, делящую интервал от 0 до 1 на $m$ промежутков.
Отсюда получаем $p_i=\tfrac i{i+j}$.
\qed

\medskip
\noindent\textit{Решение второй задачи.}
Давайте использовать ту же нумерацию позиций, как и в решении первой задачи.

Пусть $s_i$ есть среднее количество шагов робота с чебураханьем на последнем шаге, если он начинает с позиции под номером $i$.
Естественно предположить, что $s_0=s_{m}=0$ ---
ведь попадание на позиции $0$ и $m$ означает, что робот уже чебурахнулся.

Попробуем, как и раньше, посчитать значение $s_i$ новым способом.
После одного шага с $i$-ой позиции
Чебуратор окажется на позиции $i-1$ или $i+1$ с равными вероятностями $\tfrac12$.
После этого ему останется в среднем пройти $s_{i-1}$ и $s_{i+1}$ шагов соответственно. 
Не забыв учесть уже пройденный шаг, получаем
\[s_i=\tfrac12\cdot(s_{i-1}+1)+\tfrac12\cdot(s_{i+1}+1),\]
или 
\[s_{i+1}=2\cdot s_i-s_{i-1}-2.\]

Применяя эти равенства, последовательно двигаясь по позициям слева направо, получаем
$s_i=i\cdot(s_1-i+1)$
для любого $i$.
Поскольку $s_m=0$, получаем $s_1=m-1$ и, значит,
$s_i=i\cdot j$.
\qed
\medskip

Заметим, что если Чебуратор начинает на середине стола, то есть $i=j$,
то он чебурахнется в среднем за $i^2$ шагов.

\medskip

Для закрепления материала мы советуем решить следующую задачу.

\medskip
\noindent\textbf{Задача.}
Пьяница вышел из бара, расположенного в 10 кварталах от своего дома на той же улице, и идёт домой.
Чтобы сделать дорогу веселее, он разнообразит её следующим образом. 
Дойдя до перекрёстка, он подкидывает монетку --- и, если она выпадает орлом, он продолжает путь в том же направлении; если же она выпадает решкой, он разворачивается и идёт в противоположном направлении. 
Если ему случается вернуться к бару, он всегда разворачивается в сторону дома; если же он дошёл до своего дома, то завершает свой путь. 

\begin{center}
\begin{lpic}[t(1 mm),b(-1 mm),r(0 mm),l(0 mm)]{pics/bar(1)}
\lbl[b]{102,4;\home}
\lbl[b]{12,4;\man}
\lbl[b]{2,4;\beer}
\end{lpic}
\end{center}

Докажите, что в среднем за такую прогулку пьяница проходит 100 кварталов, а среднее число его возвращений к бару равно 9.

\section{Бактерии в пробирке}

\noindent\textbf{Задача.}
В пробирке живут 10 бактерий: 3 зелёных и 7 жёлтых.
Каждую минуту происходит следующее: равновероятно выбирается пара из 10 бактерий, после чего первая бактерия пары погибает, а вторая делится на две своих точных копии.
Таким образом, в конце минуты в пробирке снова остаётся ровно 10 бактерий.

Какова вероятность того, что через некоторое время все 10 бактерий будут зелёными?

\medskip
\noindent\textit{Решение.}
Мы сведем задачу к уже решённой.

Предположим, что через несколько минут 
число зелёных бактерий стало $i$ (где $0<i<10$).
Тогда спустя минуту их число может стать $i-1$, $i$ или $i+1$.
При этом вероятности первого и последнего исходов равны.
Действительно, исход с $i-1$ зелёной бактерией происходит, когда первая
бактерия случайно выбранной пары оказалась зелёной, а вторая --- жёлтой.
Исход с $i+1$ зелёной бактерией происходит, наоборот, когда первая бактерия случайно выбранной пары --- жёлтая, вторая --- зелёная. Эти две ситуации симметричны и случаются с одинаковыми вероятностями.
Исход $i$, когда выбраны бактерии одного цвета, можно считать «пропусканием хода». 

Таким образом, наша задача становится похожей на задачу про робота Чебуратора --- только теперь Чебуратор ходит направо и налево с равными положительными вероятностями, меньшими, чем $\tfrac12$, и с оставшейся вероятностью пропускает ход.
Однако заметим, что пропускание хода можно вовсе не учитывать.
То есть, ответ в задаче можно получить, подставив $i=3$ и $j=7$ в задаче про Чебуратора на длинном столе.

(Точнее, состояние, когда в пробирке ровно $i$ зелёных бактерий
соответствует $i$-ой позиции Чебуратора на длинном столе; пометим это состояние номером $i$.
Два состояния одноцветности бактерий --- состояния под номером 0 и под номером 10 ---  соответствуют состояниям робота, чебурахнувшегося с левого и правого края стола; состояния под номерами от 1 до 9 будут соответствовать всем 9-ти позициям робота на столе.)

Таким образом, с вероятностью $\tfrac{3}{10}$ все бактерии станут зелёными, 
и с вероятностью $\tfrac{7}{10}$ все станут жёлтыми.
\qed
\medskip

Эту же задачу (а значит, и первую задачу про Чебуратора)
можно решить без вычислений.

Сначала докажем, что с вероятностью 1
через некоторое время все бактерии в пробирке будут потомками одной.

Для начала заметим, что с положительной вероятностью это может случиться за первые 10 минут. 
Вероятность этого события (обозначим это число $p$) мала,
но вероятность того, что за достаточно длинный промежуток времени найдутся такие 10 минут, в пределе равна 1. 

Действительно, вероятность противоположного исхода (того, что после 10-ти минут в пробирке будут потомки более чем одной из первоначально находившихся там бактерий) равна $1-p$, и это число меньше 1. 
Значит, вероятность того, что никакой из $n$ последовательных 10-минутных интервалов не завершился тем, что в пробирке остались только потомки одной из бактерий, бывших в пробирке в начале этого интервала, равна $(1-p)^n$. 
При увеличивающемся $n$ это число приближается сколь угодно близко к нулю, значит, предельная вероятность того, что рано или поздно такой 10-минутный интервал случится --- при неограниченном увеличении времени ожидания такого интервала --- равна 1. 

При этом в самом начале каждая из десяти бактерий 
имеет равные шансы на то, чтобы стать прародителем всех выживших.
Поэтому в трех из десяти случаев все станут зелёными, 
и в семи из десяти все станут жёлтыми.

\section{Броуновское движение}

Рассмотрим случай, когда Чебуратор стоит на середине длинного стола с чётным $m=2\cdot n$ ---
чтобы чебурахнутся, ему необходимо сделать или $n$ шагов вправо, или столько же влево.
Как мы выяснили, среднее количество шагов, которое он делает, чебурахаясь на последнем из них, равно $n^2$. 
Это наблюдение можно использовать в обратном направлении:
повторив это испытание много раз и оценив среднее значение числа шагов до чебураханья, мы, 
зная размер стола, сможем оценить длину шага Чебуратора.
Последнее может оказаться полезным, если шаги Чебуратора настолько мелки, что не поддаются прямому измерению. 

Поведение Чебуратора на столе мало чем отличается 
от хаотического движения малых твёрдых частиц, плавающих в жидкости, под действием ударов молекул этой жидкости. 
Это хаотическое движение открыл Ян Ингенхауз в конце XVIII-го века,
позже его переоткрыл Роберт Броун;
в начале XX-го века удалось оценить
среднее значение числа молекул в единице объема по параметрам броуновского движения.
Этот эксперимент послужил убедительным аргументом в утверждении атомной теории строения вещества; при этом для его проведения оказалось достаточным
наблюдений с помощью обычного микроскопа.

Конечно, эта физическая задача сложней наших задач о Чебураторе.
Тем не менее, идея решения у этих задач одна и та же.

\section{Да, а о чём мы говорили?}

Мы рассмотрели частный случай \emph{цепей Маркова},
названных так в честь Андрея Андреевича Маркова-старшего (его сына тоже звали Андреем Андреевичем и он тоже занимался математикой).

А именно, мы рассматривали \emph{конечные цепи Маркова с поглощающими состояниями}.
Каждое \emph{поглощающее состояние} соответствует чебураханью нашего Чебуратора --- с левого или с правого края стола.

Цепь Маркова задаётся ориентированным графом, 
на каждом ребре которого написана вероятность перехода из одного состояния в другое.
Например, задача о Чебураторе на коротком столе может быть описана таким графом:
\begin{center}
\begin{tikzpicture}[start chain=going right]
\node[state, on chain]                 (0) {0};
\node[state, on chain]                 (1) {1};
\node[state, on chain]                 (2) {2};
\node[state, on chain] (3) {3};

\draw[
    >=latex,
    auto=right,                      
    loop above/.style={out=75,in=105,loop},
    every loop,
    ]
     (3)   edge[loop right] node {$1$}   (3)
     (2)   edge[bend left]           node[ swap] {$\tfrac12$}   (1)
           edge              node [ swap] {$\tfrac12$}   (3)
     (1)   edge             node {$\tfrac12$}   (0)
           edge  [bend left]           node [ swap]{$\tfrac12$}   (2)
     (0)   edge[loop left] node {$1$}   (0);
\end{tikzpicture}
\end{center}
Состояниям 0, 1, 2 и 3 соответствуют состояния Чебуратора «чебурахнулся налево», «стоит на левом краю стола», «стоит на правом краю стола» и «чебурахнулся направо».
Заметим, что сумма вероятностей на стрелках, исходящих из какого-то состояния, всегда равна единице.

Ключевым свойством цепи Маркова является то, что на каждом шаге
мы полностью забываем «прошлое».
Иначе говоря, «будущее» не зависит от «прошлого» при известном «настоящем».

Цепи Маркова оказались удобной и очень востребованной моделью.
Например, ранжирование веб-страниц, используемое поисковыми системами, основано на цепях Маркова.
С помощью цепей Маркова создаются алгоритмы для автоматического написания музыки.
Есть важные приложения в химии, биологии, лингвистике, экономике --- а одно из применений этой теории в физике мы вкратце описали в предыдущем разделе.

Про цепи Маркова можно почитать в популярной книжке \cite{dynkin-uspenskij}.
Затронутая нами тема \emph{случайных блужданий} рассматривается более подробно в популярных книжках и статьях 
\cite{KZhP}, 
\cite{mosteller}, 
\cite{sobolev}, 
\cite{SKU}; 
последняя статья основана на интересной физической интерпретации случайных блужданий.

\end{document}